# CURVE CROSSING FOR RANDOM WALKS REFLECTED AT THEIR MAXIMUM[1]

### By Ron Doney and Ross Maller

*University of Manchester and Australian National University*

Let $R_n = \max_{0 \le j \le n} S_j - S_n$ be a random walk $S_n$ reflected in its maximum. Except in the trivial case when $P(X \ge 0) = 1$, $R_n$ will pass over a horizontal boundary of any height in a finite time, with probability 1. We extend this by giving necessary and sufficient conditions for finiteness of passage times of $R_n$ above certain curved (power law) boundaries, as well. The intuition that a degree of heaviness of the negative tail of the distribution of the increments of $S_n$ is necessary for passage of $R_n$ above a high level is correct in most, but not all, cases, as we show. Conditions are also given for the finiteness of the expected passage time of $R_n$ above linear and square root boundaries.

**1. Introduction and preliminary results.** Let $X, X_1, X_2, \ldots,$ be i.i.d. r.v.'s, not degenerate at 0, with c.d.f. $F(\cdot)$ on $\mathbb{R}$, and

$$S_n = X_1 + X_2 + \cdots + X_n, \qquad S_0 = 0,$$

the corresponding random walk. Denote by

$$R_n = \max_{0 \le j \le n} S_j - S_n, \qquad n = 0, 1, 2, \ldots,$$

the random walk reflected in its maximum. Of course, $R_n \ge 0$, $n = 0, 1, 2, \ldots$.

The reflected process $R_n$ is of fundamental importance in the theory of random walks and is also an object of interest, in itself, in many applied areas, such as queueing theory; see, for example, [1, 12, 18, 29] and their references. More recently, $R_n$ has been used extensively in various other kinds of modeling. The first time the reflected process upcrosses a fixed level gives the optimal time to exercise a "Russian" option [2, 26, 28]. Hansen [16]

---

Received September 2005; revised October 2006.

[1]Supported in part by ARC Grants DP0210572 and DP0664603.

---

has some interesting generalizations and an application to genetics of the maximal sequence $R_n^* := \max_{1 \le j \le n} R_j$. There are many other applications of $R_n$ in finance studies and elsewhere. See also [10, 22] and [23].

$R_n$ has been intensively studied in conjunction with these applications, but its renewal-theoretic properties *per se* seem to have received little attention so far. Here we consider some very basic but important questions related to this aspect. Thus, in Theorem 2.1 of Section 2, we give necessary and sufficient conditions for the almost sure (a.s.) finiteness of passage times of $R_n$ out of power law regions of the form $\{(n, y) : 0 \le y \le rn^\kappa, n = 1, 2, \ldots\}$, where $r \ge 0$, $\kappa = 0$, or $r > 0$, $\kappa > 0$. Then, in Theorem 2.2, we give conditions for the finiteness of expected values of passage times of $R_n$ out of linear ($\kappa = 1$) or parabolic ($\kappa = 1/2$) regions. These can be thought of as extensions or generalizations of similar results for random walks, and we use a variety of the techniques developed for random walks in their proofs.

To complete the present section, we introduce some notation which will be useful throughout the paper, and state an introductory Proposition 1.1 which helps to motivate the kinds of issues we will consider. Let

$$(1.1) \qquad S_n^* = \max_{0 \le j \le n} S_j, \qquad n = 0, 1, 2, \ldots.$$

Then it is easy to see that

$$(1.2) \qquad R_n = S_n^* - S_n = \max_{0 \le j \le n} (S_j - S_n) \stackrel{D}{=} - \min_{0 \le j \le n} S_j.$$

The identity (1.2) (equality in distribution for each $n = 1, 2, \ldots$, but not of processes) is of course well known. Another useful representation is to write $R_n$ as the sum of its increments:

$$(1.3) \qquad R_n = \sum_{i=1}^n \Delta_i,$$

where, as is easily checked,

$$(1.4) \quad \Delta_i = R_i - R_{i-1} = -X_i \mathbf{1}_{\{X_i \le R_{i-1}\}} - R_{i-1} \mathbf{1}_{\{X_i > R_{i-1}\}}, \qquad i = 1, 2, \ldots.$$

Note that, if $F(0-) = 0$, then $R_n$ is identically 0, while if $F(0) = 1$, then $R_n = -S_n$, the negative of a random walk. The first case is trivial and, for the second, the results we examine are already known, as discussed later (and, in fact, our present results remain true in this case, with appropriate interpretations), so we exclude them in what follows. Thus, *throughout, we make the blanket assumption that $0 < F(0-) \le F(0) < 1$.* Throughout, also, we will use "r.v." to mean "random variable," "$\stackrel{D}{\to}$" for convergence in distribution, "$\stackrel{P}{\to}$" for convergence in probability, "a.s." for almost sure convergence and "i.o." for "infinitely often." Let $X^+ = \max(X, 0)$ and $X^- = X^+ - X$ (and similarly for $X_i^+$ and $X_i^-$).

Our first proposition lists some elementary properties of $R_n$.

PROPOSITION 1.1.    (a) $\liminf_{n\to\infty} P(R_n > x) > 0$ for every $x > 0$; consequently, we always have $\limsup_{n\to\infty} R_n = +\infty$ a.s., and we never have $R_n \overset{P}{\to} 0$ ($n \to \infty$).

(b) (i) $R_n \overset{P}{\to} \infty$ ($n \to \infty$) if and only if $\liminf_{n\to\infty} S_n = -\infty$ a.s.

(ii) The following four conditions are equivalent: $P(R_n = 0 \ i.o.) < 1$; $\lim_{n\to\infty} S_n = -\infty$ a.s.; $\lim_{n\to\infty} R_n = +\infty$ a.s.; and

$$(1.5) \qquad \sum_{n\geq 1} P(R_n \leq x) < \infty \qquad \text{for some (hence, every) } x \geq 0.$$

REMARKS.    (i) It is easy to show from (1.2) that $R_n$ is tight as $n \to \infty$ if and only if $\lim_{n\to\infty} S_n = +\infty$ a.s., and, in fact, this implies that $R_n \overset{D}{\to} R$, where $R \overset{D}{=} -\min_{0 \leq j < \infty} S_j$, with $P(0 \leq R < \infty) = 1$. Thus, $R_n$ is stochastically bounded if and only if $S_n$ drifts to $+\infty$ a.s. This situation has been well studied in various applications (see, e.g., [3], page 388 and [29]), and we will mainly be concerned with the other cases, when $S_n$ oscillates or drifts to $-\infty$ a.s., so that $R_n$ continues to grow with $n$ [part (b) of Proposition 1.1].

(ii) Analytic conditions for $\lim_{n\to\infty} S_n = \pm\infty$ a.s. are in [20]. See Proposition 2.1 below for $\liminf_{n\to\infty} S_n = -\infty$ a.s.

(iii) We remark that, with the obvious modifications, all of our results apply to the reflected process $r_n := S_n - \min_{0 \leq j \leq n} S_j$. For a financial application of $r_n$, see [13].

## 2. Passage times above power law boundaries.

We can measure the rate of growth of $R_n$ by seeing how quickly it leaves a region. We restrict ourselves here to power law regions. Thus, for constants $\kappa > 0$, $r > 0$, or $\kappa = 0$, $r \geq 0$, define

$$(2.1) \qquad \tau_\kappa(r) = \min\{n \geq 1 : R_n > rn^\kappa\}.$$

(Throughout, give the minimum of the empty set the value $+\infty$.) Basic questions of interest are to find conditions on $F$ which are equivalent to $\tau_\kappa(r)$ being a.s. finite or having a finite expectation. For random walk, the first question is answered in [21] and [9]; a summary of their results (with a sign change) is in Proposition 2.1, later in this section. We build on these to give our first main result for $R_n$. It might seem obvious, a priori, in keeping with Proposition 1.1, that a certain heaviness of the *negative* tail of $F$ is required in order for $R_n$ to escape the power law region. However, when $\kappa \in (1/2, 1)$, just as in the case of a random walk crossing a one-sided boundary, this intuition can fail. The second part of (2.3) below can hold even when $X$ is stochastically bounded below, so that the negative tail of $F$ is zero for large $x$; see part (e)(ii) of Proposition 2.1.

Thus, delineating the precise conditions is not at all straightforward. We find the following:

THEOREM 2.1. (a) *Suppose* $\kappa = 0$. *We have* $\tau_0(r) < \infty$ *a.s. for all* $r \geq 0$, *and, in fact,* $E(e^{\lambda \tau_0(r)}) < \infty$, *at least for some small enough* $\lambda > 0$, *for all* $r \geq 0$.

(b) *Suppose* $\kappa > 0$. *We have* $\tau_\kappa(r) < \infty$ *a.s. for all* $r > 0$ *if and only if:*

(i) *for* $\kappa > 1$,

$$(2.2) \qquad\qquad\qquad E(X^-)^{1/\kappa} = \infty;$$

(ii) *for* $0 < \kappa \leq 1$,

$$(2.3) \qquad E(X^-)^{1/\kappa} = \infty \quad or \quad \liminf_{n \to \infty} \left( \frac{S_n}{n^\kappa} \right) = -\infty \qquad a.s.$$

Explicit criteria in terms of the distribution function $F$ of the $X_i$ for $\liminf_{n \to \infty} S_n / n^\kappa = -\infty$ a.s. are listed in Proposition 2.1 below. Parts (a) and (b) of the proposition are essentially due to [6] and [11], respectively; parts (c) and (d) are in [21]; part (e) and the following comment is from [9]. (Actually, these papers deal with the condition $\limsup_{n \to \infty} S_n / n^\kappa = +\infty$ a.s., but the results for $\liminf_{n \to \infty} S_n / n^\kappa$ follow after a sign reversal.) To state them, let $\bar{F}(y) := 1 - F(y)$, and define the integrals

$$(2.4) \qquad A_+(x) = \int_0^x \bar{F}(y)\,dy, \qquad x > 0, \quad \text{and} \quad J_- = \int_{[0,\infty)} \frac{x|dF(-x)|}{A_+(x)}.$$

Note that $0 \leq A_+(x) \leq EX^+$. We let $A_+(x)/x$ have its limiting value, $1 - F(0) > 0$, at 0. We also need the function defined, for $y \geq 0$, when $EX^+ < \infty$, as

$$(2.5) \quad W(y) = \int_0^y \int_{(z,\infty)} xF(dx)\,dz = y \int_{(y,\infty)} zF(dz) + \int_{[0,y]} z^2 F(dz).$$

Note that $W(y) > 0$ for all $y > 0$, since we assume that $F$ is not concentrated on $(\infty, 0]$. Define, for $\lambda > 0$, $y > 0$, and $1/2 < \kappa < 1$,

$$(2.6) \qquad I_\kappa(\lambda) := \int_1^\infty \exp\left\{ -\lambda \left( \frac{y^{(2\kappa-1)/\kappa}}{W(y)} \right)^{\kappa/(1-\kappa)} \right\} \frac{dy}{y} \leq \infty,$$

and let

$$(2.7) \qquad\qquad \lambda_\kappa^* = \inf\{\lambda > 0 : I_\kappa(\lambda) < \infty\} \in [0, \infty].$$

PROPOSITION 2.1. $\liminf_{n \to \infty} S_n / n^\kappa = -\infty$ *a.s. if and only if:*

(a) *when* $\kappa > 1$,

$$(2.8) \qquad\qquad \int_{[1,\infty)} \left( \frac{x^{1/\kappa}}{1 + x^{(1/\kappa)-1} A_+(x)} \right) |dF(-x)| = \infty;$$

(b) *when $\kappa = 1$ or $\frac{1}{2} < \kappa < 1$ and $E|X| = \infty$,*

$$(2.9) \qquad\qquad J_- = \infty;$$

(c) *when $0 \leq \kappa \leq \frac{1}{2}$,*

$$(2.10) \qquad\qquad J_- = \infty \quad or \quad 0 \leq -EX \leq E|X| < \infty;$$

(d) *when $\frac{1}{2} < \kappa < 1$, $E|X| < \infty$ and $EX \neq 0$,*

$$(2.11) \qquad\qquad EX < 0;$$

(e) *when $\frac{1}{2} < \kappa < 1$, $E|X| < \infty$ and $EX = 0$,*

$$(2.12) \qquad\qquad (i) \ E(X^-)^{1/\kappa} = \infty, \quad or$$

$$(2.13) \qquad\qquad (ii) \ E(X^-)^{1/\kappa} < \infty = E(X^+)^{1/\kappa} \ and \ \lambda_\kappa^* = \infty.$$

*Furthermore, when $\frac{1}{2} < \kappa < 1$, $E|X| < \infty$, $EX = 0$ and $E(X^-)^{1/\kappa} < \infty = E(X^+)^{1/\kappa}$, then $\liminf_{n\to\infty} S_n/n^\kappa = 0$ a.s. if and only if $\lambda_\kappa^* = 0$.*

REMARKS.  (i) Our blanket assumption that $0 < F(0-) \leq F(0) < 1$ is not restrictive in Proposition 2.1, because if $F(0) = 1$, then $\liminf_{n\to\infty} S_n/n^\kappa = -\limsup_{n\to\infty} |S_n|/n^\kappa$ a.s. and Theorem 1 of [21] gives the required criterion; while if $F(0-) = 0$, then neither $\liminf_{n\to\infty} S_n/n^\kappa = -\infty$ a.s. nor any of (2.8)–(2.13) can occur.

(ii) In general, neither of the two conditions in (2.3) imply each other, as can be seen from a perusal of Proposition 2.1.

(iii) Again, given our assumption that $0 < F(0-) \leq F(0) < 1$, the a.s. finiteness of $\tau_\kappa(r) < \infty$ a.s. is equivalent to $\limsup_{n\to\infty} R_n/n^\kappa > r$ a.s. (see Lemma 3.1 of Section 3). Thus, the contrapositives of the conditions in Theorem 2.1 give equivalences for $\limsup_{n\to\infty} R_n/n^\kappa$ to be finite a.s. We summarize these in the following:

COROLLARY 2.1 (Corollary to Theorem 2.1).  *$\limsup_{n\to\infty} R_n/n^\kappa$ is finite a.s. if and only if:*

(a) *when $\kappa \geq 1$, $E(X^-)^{1/\kappa} < \infty$;*
(b) *when $0 < \kappa \leq 1$, $E(X^-)^{1/\kappa} < \infty$ and $\liminf_{n\to\infty} S_n/n^\kappa > -\infty$ a.s.*

In the case $\kappa = 0$, we have $\limsup_{n\to\infty} R_n = \infty$ a.s., by Proposition 1.1, because we always assume that $F(0-) > 0$. From the proof of Theorem 2.1, it can further be seen that $\limsup_{n\to\infty} R_n/n^\kappa$, when finite a.s., is in fact 0 a.s., except in the cases (1) $\kappa = 1$, $E|X| < \infty$, $EX < 0$ (when $\lim_{n\to\infty} R_n/n = |EX|$ a.s. by the strong law of large numbers for $S_n$) and (2) $1/2 < \kappa < 1$, $E|X| < \infty$, $EX = 0$ and $E(X^-)^{1/\kappa} < \infty = E(X^+)^{1/\kappa}$, when $\liminf_{n\to\infty} R_n/n^\kappa \in (0, \infty)$ a.s. if and only if $\lambda_\kappa^* \in (0, \infty)$; see the remark following the statement of Lemma 3.4.

Our second main result considers the expected value of the passage time of $R_n$ above linear and square root boundaries.

THEOREM 2.2.  (a) *Suppose* $\sigma^2 = EX^2 < \infty$ *and* $EX = 0$. *Then:*

  (i) $E\tau_{1/2}(\sigma r) = \infty$ *for* $r \geq 1$;
  (ii) $E\tau_{1/2}(\sigma r) < \infty$ *for* $r < 1$.

 (b) *Suppose* $E|X| < \infty$ *and* $EX < 0$. *Then:*
    (i) $E\tau_1(r) < \infty$ *for* $r < -EX$;
    (ii) *If* $r > -EX$, $P(\tau_1(r) = \infty) > 0$.

 (c) *Suppose* $E|X| < \infty$ *and* $EX < 0$, *and, in addition,* $E(X^+)^2 < \infty$. *Then* $E\tau_1(r) = \infty$ *for* $r = -EX$.

The random walk precursor of Theorem 2.2 (a) is in [4] and [15], who dealt with independent orthonormal r.v.'s (having mean 0 and finite variance), and showed that the first time ($T_r$, say) at which the corresponding partial sum exits the parabolic region $\{(n, y) : |y| \leq r\sqrt{n}, n = 1, 2, \ldots\}$, where $r > 0$, has finite expectation if $r < EX^2$ and infinite expectation otherwise.

A corresponding linear version is the following: suppose $E|X| < \infty$ and $EX > 0$, then the first time a random walk with step $X$ starting from 0 passages above the line $y = rn$, $n = 1, 2, \ldots$, $r \geq 0$, has finite expectation if and only if $r < EX$. This is easily proved by reducing the problem to the finiteness or otherwise of the expected first passage time above 0 of a random walk which is drift free or has negative drift when considering the cases $r \geq EX$ and has positive drift otherwise; see, for example, [14] for a discussion of this.

Some other results are not so easily settled, even in the random walk case. For example, Gundy and Siegmund [15] conjecture that, in the above notation, $ET_r = \infty$ continues to hold when $EX^2 = \infty$, for all $r > 0$. They mention having a proof of this for the case when the $X_i$ are symmetrically distributed, but the general problem remains open. Likewise, in our Theorem 2.2, the restriction $E(X^+)^2 < \infty$ may not be necessary in part (c).

A natural extension of our results is to ask for conditions for the finiteness or otherwise of $E\tau_\kappa(r)$ when $\kappa \neq 1/2$ or 1. Again, in view of the above discussion, we expect this may be a rather difficult exercise. But $\kappa = 1/2$ or 1 are probably the most important practical cases.

We refer to Novikov [24, 25] and his references for more precise estimates of magnitudes of tail probabilities of stopping time distributions, under certain conditions.

CONCLUDING REMARKS.  As might be expected, there is a counterpart of Theorem 2.1 relating to the large time behavior of a Lévy process, and also for the results of Proposition 1.1, with appropriate interpretations. The proofs can be constructed as in [8, 9], using the methods of [7]. We omit the details. Lévy versions of Theorem 2.2 have been proved by Savov [27].

**3. Proofs.** Recall our blanket assumption throughout that $0 < F(0-) \leq F(0) < 1$.

PROOF OF PROPOSITION 1.1. (a) Since $F(0-) > 0$, there are $\varepsilon > 0$, $\delta > 0$, such that $P(X \leq -\varepsilon) > \delta$. Take any $x > 0$ and choose $K > 1$ so that $K\varepsilon > x$. Suppose $X_i \leq -\varepsilon$, $n+1 \leq i \leq n+K$. Then for $n = 1, 2, \ldots$,

$$S_{n+K} = S_n + \sum_{i=n+1}^{n+K} X_i \leq S_n - K\varepsilon < S_n - x \leq S_{n+K}^* - x,$$

so $R_{n+K} > x$. Thus, for $n = 1, 2, \ldots$,

$$P(R_{n+K} > x) \geq P(X_i \leq -\varepsilon, n+1 \leq i \leq n+K) \geq \delta^K > 0,$$

so $\liminf_{n \to \infty} P(R_n > x) > 0$. It then follows from the Hewitt–Savage 0–1 law ([17] or [5], page 226) that $\limsup_{n \to \infty} R_n = +\infty$ a.s., and clearly, also, $R_n \overset{P}{\to} 0$ $(n \to \infty)$ cannot occur.

(b) (i) We have

$$
\begin{aligned}
R_n \overset{P}{\to} \infty \quad &\Longleftrightarrow \quad \lim_{n \to \infty} P(R_n \leq x) = 0 \qquad \text{for all } x > 0 \\
&\Longleftrightarrow \quad \lim_{n \to \infty} P\Big(\min_{0 \leq j \leq n} S_j \geq -x\Big) = 0 \qquad [\text{by } (1.2)] \\
&\Longleftrightarrow \quad \min_{0 \leq j \leq n} S_j \overset{P}{\to} -\infty \\
&\Longleftrightarrow \quad \min_{0 \leq j \leq n} S_j \to -\infty \qquad \text{a.s. (since the sequence is monotone)} \\
&\Longleftrightarrow \quad \liminf_{n \to \infty} S_n = -\infty \qquad \text{a.s.}
\end{aligned}
$$

(ii) Let $P(R_n = 0 \text{ i.o.}) < 1$ and suppose $S_n$ does not drift to $-\infty$ a.s. Then $\limsup_{n \to \infty} S_n = +\infty$ a.s. and so there are infinitely many ascending ladder times, a.s. This means that $S_n$ exceeds $S_{n-1}^*$ infinitely often, a.s., hence, $S_n^* = S_n$ i.o. a.s., and so $R_n = 0$ i.o. a.s., a contradiction. Thus, $\lim_{n \to \infty} S_n = -\infty$ a.s. Next, $\lim_{n \to \infty} S_n = -\infty$ a.s. implies $R_n = S_n^* - S_n \geq -S_n \to \infty$ a.s., while $R_n \to \infty$ a.s. obviously implies $P(R_n = 0 \text{ i.o.}) = 0$.

For the final equivalence, assume $\sum_n P(R_n \leq x) < \infty$ for some $x \geq 0$. Then $\sum_n P(R_n = 0) < \infty$, so by the Borel–Cantelli lemma, $P(R_n = 0 \text{ i.o.}) = 0$. This implies $\lim_{n \to \infty} S_n = -\infty$ a.s., as just shown, and this further implies, by Theorem 2.1 of [20] [interchanging $+$ and $-$ in their result, i.e., applying their result to the random walk $\widetilde{S}_n = \sum_{i=1}^n (-X_i)$], that

$$\infty > \sum_{n \geq 1} P\Big(\min_{0 \leq j \leq n} S_j \geq -x\Big) = \sum_{n \geq 1} P(R_n \leq x) \qquad \text{for every } x \geq 0. \qquad \square$$

PROOF OF THEOREM 2.1. (a) Take an $r \geq 0$. As in the proof of Proposition 1.1, since $F(0-) > 0$, we can choose $\varepsilon > 0$, $\delta > 0$, $K \geq 1$, so that $F(-\varepsilon-) \geq \delta$ and $K\varepsilon > r$. Then for $n = 0, 1, \ldots,$

$$P(\tau_0(r) \leq n + K | \tau_0(r) > n) \geq P(X_i \leq -\varepsilon, \ n+1 \leq i \leq n+K) \geq \delta^K,$$

and the required results in part (a) both follow from

$$P(\tau_0(r) > nK) = P(\tau_0(r) > K) \prod_{j=2}^{n} P(\tau_0(r) > jK | \tau_0(r) > (j-1)K)$$

$$\leq (1 - \delta^K)^n.$$

The following lemma will be useful in the rest of the proof.

LEMMA 3.1. *Take $r > 0$ and $\kappa > 0$. Then $\limsup_{n \to \infty} R_n / n^\kappa > r$ a.s. if and only if $\tau_\kappa(r) < \infty$ a.s.*

PROOF. Since $\{\tau_\kappa(r) > n\} = \{\max_{1 \leq j \leq n}(R_j / j^\kappa) \leq r\}$, we have $P(\tau_\kappa(r) = \infty) = P(\max_{j \geq 1}(R_j / j^\kappa) \leq r)$. Now $\limsup_{n \to \infty} R_n / n^\kappa > r$ a.s. implies $P(\max_{j \geq 1}(R_j / j^\kappa) > r) = 1$, so one direction of the proof is obvious.

Conversely, suppose $\tau_\kappa(r) < \infty$ a.s., so $\lim_n P(\max_{1 \leq j \leq n}(R_j / j^\kappa) \leq r) = 0$. Note that then

$$P\left(\max_{1 \leq j \leq n}(R_j / j^\kappa) \leq r \text{ i.o.}\right) = \lim_n P\left(\max_{1 \leq j \leq m}(R_j / j^\kappa) \leq r \text{ for some } m > n\right)$$

$$\leq \lim_n P\left(\max_{1 \leq j \leq n}(R_j / j^\kappa) \leq r\right) = 0.$$

We wish to show $P(\max_{k \leq j \leq n}(R_j / j^\kappa) \leq r \text{ i.o.}) = 0$ for each $k \geq 1$, and proceed by induction. Let $A_n^k := \{\max_{k \leq j \leq n}(R_j / j^\kappa) \leq r\}$. Suppose $P(A_n^k \text{ i.o.}) = 0$ for some $k \geq 1$, but $P(A_n^{k+1} \text{ i.o.}) > 0$. Then by the Hewitt–Savage law,

$$1 = P(A_n^{k+1} \text{ i.o.})$$

$$= \lim_m P(A_n^{k+1}, \text{ for some } n > m)$$

$$\leq \lim_m P(A_n^{k+1}, \text{ for some } n > m, R_k \leq rk^\kappa) + P(R_k > rk^\kappa)$$

$$\leq \lim_m P(A_n^k, \text{ for some } n > m) + P(R_k > rk^\kappa)$$

$$= P(A_n^k \text{ i.o.}) + P(R_k > rk^\kappa)$$

$$= P(R_k > rk^\kappa),$$

giving $R_k > rk^\kappa$ a.s. This is not possible when $F(0-) > 0$ since then $P(R_1 = 0) = P(X \geq 0) > 0$, thus, $P(R_k = 0) > 0$ for each $k \geq 1$. So $P(\max_{k \leq j \leq n}(R_j / j^\kappa) \leq r \text{ i.o.}) = 0$, $k \geq 1$, thus, $\lim_n P(\max_{k \leq j \leq n}(R_j / j^\kappa) \leq r) = P(\max_{j \geq k}(R_j / j^\kappa) \leq$ $r) = 0$, $k \geq 1$. Letting $k$ tend to $\infty$ then gives $P(R_j/j^\kappa > r$ i.o.$) = 1$, which implies $\limsup_{n\to\infty} R_n/n^\kappa > r$ a.s.  $\square$

We now return to the proof of Theorem 2.1. We first prove the forward direction for both parts of (b).

(b) (i) Keep $\kappa > 1$, and suppose $\tau_\kappa(r) < \infty$ a.s. for some $r > 0$. If $E(X^-)^{1/\kappa} < \infty$, the Marcinkiewicz–Zygmund law (e.g., [5], page 125) gives $\lim_{n\to\infty}(\sum_{i=1}^n X_i^-/n^\kappa) = 0$, a.s. If this is so, then

$$R_n = \max_{0 \leq j \leq n} S_j - S_n = \max_{0 \leq j < n}\left(-\sum_{i=j+1}^n X_i\right) \vee 0$$

$$\leq \max_{0 \leq j < n}\left(\sum_{i=j+1}^n X_i^-\right) = \sum_{i=1}^n X_i^- = o(n^\kappa) \qquad \text{a.s.}$$

But by Lemma 3.1, this is a contradiction. Thus, the forward direction of part (i) is proved.

(ii) Keep $0 < \kappa \leq 1$. Let $T_n$ be the strict increasing ladder times of $S_n$, that is, $T_0 = 0$ and

$$(3.1) \qquad T_n = \min\{j > T_{n-1} : S_j > S_{T_{n-1}}\}, \qquad n = 1, 2, \dots.$$

If $T_{n-1} < \infty$, define the *depth* of an excursion of $S_n$ below the maximum as

$$(3.2) \qquad D_n = \max_{T_{n-1} \leq j < T_n}\left(-\sum_{i=T_{n-1}+1}^j X_i\right), \qquad n = 1, 2, \dots.$$

[In (3.2), and throughout, we make the convention that $\sum_{i=a}^b = 0$ when $b < a$.] The r.v. $D_n$ measures the *height* of an excursion of $R_n$ away from 0; we have $R_{T_n} = 0$, $n = 1, 2$ and

$$(3.3) \qquad \max_{T_{n-1} \leq j < T_n} R_j = D_n, \qquad n = 1, 2, \dots.$$

[If two ladder times $T_{n-1}, T_n$ occur at consecutive integers, so that $R_{T_{n-1}} = R_{T_n} = 0$, (3.2) gives $D_n = 0$, agreeing with (3.3), and formally registering that the depth of the nonexistent excursion is 0.]

LEMMA 3.2.  *Keep $0 < \kappa \leq 1$ and suppose $\lim_{n\to\infty} S_n = +\infty$ a.s. Then $E(X^-)^{1/\kappa} < \infty$ if and only if $E(D_1^{1/\kappa}) < \infty$.*

PROOF.  Assume $\lim_{n\to\infty} S_n = +\infty$ a.s. Then $T_n < \infty$ a.s. for all $n$, and, in fact, $ET_1 < \infty$; see, for example, Theorem II.9.1, page 66, in [14]. Thus, the $D_n$ are well defined. Since $S_j \leq 0, 0 \leq j < T_1$, we have

$$D_1 = \max_{0 \leq j < T_1}(-S_j) \geq -S_1 = S_1^- = X_1^-,$$

and one direction of the lemma is obvious. Conversely,

$$D_1 \leq \max_{1 \leq j < T_1} \left( \sum_{i=1}^{j} X_i^- \right) = \sum_{i=1}^{T_1-1} X_i^- = F_1 \qquad \text{say.}$$

Now for $0 < \kappa \leq 1$, $E(X^-)^{1/\kappa} < \infty$ and $\lim_{n \to \infty} S_n = +\infty$ a.s. imply $ET_1^{1/\kappa} < \infty$ ([20], Theorem 2.1), so we can apply Theorem I.5.2, page 22, in [14] to get $EF_1^{1/\kappa} < \infty$ and, hence, $ED_1^{1/\kappa} < \infty$. $\square$

LEMMA 3.3. *Keep $\kappa > 0$. If $ED_1^{1/\kappa} < \infty$ and $\lim_{n \to \infty} S_n = +\infty$ a.s., then $\lim_{n \to \infty} R_n / n^\kappa = 0$ a.s., and so $P(\tau_\kappa(r) = \infty) > 0$ for all $r > 0$.*

PROOF. Again with $T_n$ as the strict increasing ladder times of $S_n$,

$$\max_{j \geq T_n} \left( \frac{R_j}{j^\kappa} \right) = \max_{m > n} \max_{T_{m-1} \leq j < T_m} \left( \frac{R_j}{j^\kappa} \right)$$

$$\leq \max_{m > n} \left( \frac{\max_{T_{m-1} \leq j < T_m} R_j}{T_{m-1}^\kappa} \right) = \max_{m > n} \left( \frac{D_m}{T_{m-1}^\kappa} \right). \tag{3.4}$$

Since $\lim_{n \to \infty} S_n = +\infty$ a.s., we have $ET_1 < \infty$, and thus, $\lim_{n \to \infty} T_m / m = ET_1$ a.s. is finite a.s. The $D_m$ are i.i.d., and with $ED_1^{1/\kappa} < \infty$, by hypothesis, so we have $\lim_{m \to \infty} D_m / m^\kappa = 0$ a.s. Thus, the right-hand side of (3.4) tends to 0 a.s. as $n \to \infty$, giving $\lim_{n \to \infty} R_n / n^\kappa = 0$ a.s. Then $P(\tau_\kappa(r) = \infty) > 0$ for all $r > 0$ follows from Lemma 3.1. $\square$

We can now complete the proof of the forward direction of part (b)(ii) of Theorem 2.1. We have $0 < \kappa \leq 1$ and $\tau_\kappa(r) < \infty$ for all $r > 0$, and must prove that (2.3) holds.

If $E(X^-)^{1/\kappa} = \infty$, then (2.3) holds, so suppose $E(X^-)^{1/\kappa} < \infty$. Then by Lemmas 3.2 and 3.3, we cannot have $\lim_{n \to \infty} S_n = +\infty$ a.s., consequently, $\liminf_{n \to \infty} S_n = -\infty$ a.s. So, using Proposition 2.1(c) with $\kappa = 0$, we see that (2.10) holds.

First suppose $0 < \kappa \leq 1/2$. Then by (2.10) again, we have

$$\liminf_{n \to \infty} \left( \frac{S_n}{n^\kappa} \right) = -\infty \qquad \text{a.s.,} \tag{3.5}$$

so (2.3) holds.

Next consider $1/2 < \kappa \leq 1$. We still have (2.10). If $E|X| = \infty$, then $J_- = \infty$ by (2.10), and then (3.5) holds by (2.9). If $\kappa = 1$, we can finish here because $E|X| < \infty$ cannot occur. If it did, we would have, a.s. as $n \to \infty$,

$$\frac{R_n}{n} = \frac{S_n^*}{n} - \frac{S_n}{n} \to (EX)^-.$$

Thus, if $EX = 0$, then $P(\tau_1(r) = \infty) > 0$ for all $r > 0$ by Lemma 3.1, while if $EX < 0$, then $P(\tau_1(r) = \infty) > 0$ for all $r > |EX|$, again by Lemma 3.1. Either is a contradiction.

Finally, consider $1/2 < \kappa < 1$ and $E|X| < \infty$. Then $EX \le 0$ by (2.10). If $EX < 0$, then $\lim_{n \to \infty} S_n/n = EX < 0$ a.s., so (3.5) and, hence, (2.3) holds. It remains to consider the case $EX = 0$.

The next lemma allows us to deal with this. Recall the definitions of $W(y)$, $I_\kappa(\lambda)$ and $\lambda_\kappa^*$, in (2.5), (2.6) and (2.7), respectively.

LEMMA 3.4. *Keep* $1/2 < \kappa < 1$. *Suppose that* $E|X| < \infty$, $EX = 0$, $E(X^-)^{1/\kappa} < \infty = E(X^+)^{1/\kappa}$, *and* $\lambda_\kappa^* < \infty$. *Then*

$$(3.6) \qquad \limsup_{n \to \infty} \left( \frac{R_n}{n^\kappa} \right) \le r_0 := 9 \cdot 2^\kappa (6\lambda_\kappa^* 2^\kappa)^{1-\kappa} \qquad a.s.,$$

*and consequently,* $P(\tau_\kappa(r) = \infty) > 0$ *for all* $r \ge r_0$.

REMARK. Suppose $1/2 < \kappa < 1$, $E|X| < \infty$, $EX = 0$ and $E(X^-)^{1/\kappa} < \infty = E(X^+)^{1/\kappa}$. If $\lambda_\kappa^* = \infty$, then Proposition 2.1(e)(ii), together with the fact that $R_n > -S_n$, $n = 1, 2, \ldots$, gives $\limsup_{n \to \infty} R_n/n^\kappa = \infty$ a.s., a partial converse to (3.6).

It is possible to have $I_\kappa(\lambda) = \infty$ for some but not all $\lambda > 0$, as shown in [9]. If this happens, then $\lambda_\kappa^* \in (0, \infty)$ and so $\liminf_{n \to \infty} S_n/n^\kappa < 0$ a.s., by part (f) of Proposition 2.1, and hence, we have $\limsup_{n \to \infty} R_n/n^\kappa > 0$ a.s., as well as $\limsup_{n \to \infty} R_n/n^\kappa < \infty$ a.s. So it is possible to have $\limsup_{n \to \infty} R_n/n^\kappa \in (0, \infty)$ a.s. in this case. Lemma 3.4 should be compared with Corollary 1.1 of [9].

PROOF OF LEMMA 3.4. Fix $1/2 < \kappa < 1$, suppose $E|X| < \infty$, $EX = 0$, $E(X^-)^{1/\kappa} < \infty = E(X^+)^{1/\kappa}$, and $\lambda_\kappa^* < \infty$. Then there is a $\lambda > \lambda_\kappa^*$ with $I_\kappa(\lambda) < \infty$. We keep this $\lambda$ fixed through the proof, then at the end let $\lambda \downarrow \lambda_\kappa^*$ to get (3.6).

We can write

$$
\begin{aligned}
(3.7) \qquad R_n &= \max_{0 \le j \le n} \left( -\sum_{i=j+1}^{n} X_i \right) \\
&= \max_{0 \le j \le n} \left( \sum_{i=j+1}^{n} X_i^- - \sum_{i=j+1}^{n} X_i^+ \right).
\end{aligned}
$$

For $D > 0$, we have

$$
\begin{aligned}
\sum_{i=j+1}^{n} X_i^+ &\ge D \sum_{i=j+1}^{n} \mathbf{1}_{\{X_i^+ > D\}} + \sum_{i=j+1}^{n} X_i^+ \mathbf{1}_{\{X_i^+ \le D\}} \\
&= (n-j)D - D \sum_{i=j+1}^{n} \mathbf{1}_{\{X_i^+ \le D\}} + \sum_{i=j+1}^{n} X_i^+ \mathbf{1}_{\{X_i^+ \le D\}}.
\end{aligned}
$$

Note that $EX^- = EX^+$ in our case, and recall that $\lim_{y\to\infty} y\bar{F}(y) = 0$ when $E|X| < \infty$. We will also use the function $\nu_+(x) := \int_{[0,x]} y\,dF(y)$, for $x > 0$. Some algebra then shows that

$$
\begin{aligned}
&\sum_{i=j+1}^{n} X_i^- - \sum_{i=j+1}^{n} X_i^+ \\
&\quad \leq \sum_{i=j+1}^{n} (X_i^- - EX^-) + D \sum_{i=j+1}^{n} \left( \mathbf{1}_{\{X_i^+ \leq D\}} - F(D) \right) \\
&\qquad - \sum_{i=j+1}^{n} (X_i^+ \mathbf{1}_{\{X_i^+ \leq D\}} - \nu_+(D)) + n \int_D^\infty \bar{F}(y)\,dy.
\end{aligned}
\tag{3.8}
$$

We will choose $D$ as follows. We have $W(x) > 0$ for all $x > 0$ (because $EX = 0$), $\lim_{x\to\infty} W(x)/x = 0$ (because $EX^+ < \infty$), and $\lim_{x\downarrow 0} W(x)/x = EX^+$. So, given $\delta > 0$ and $x > x_0 := (\delta/EX^+)^{1/(1-\kappa)}$, we can define $D(x) = D(x, \delta)$ by

$$
D(x) = \inf\left\{ y > 0 : \frac{W(y)}{y} \leq \frac{\delta}{x^{1-\kappa}} \right\}.
$$

Then $0 < D(x) < \infty$ for $x > x_0$, $\lim_{x\to\infty} D(x) = \infty$, and by the continuity of $W(x)$, $D(x)$ satisfies

$$
\frac{x^{1-\kappa} W(D(x))}{D(x)} = \delta.
\tag{3.9}
$$

Now take $k \geq 1$ and $1 \leq n \leq 2^k$, and let

$$
A_n = \sum_{i=1}^{n} (X_i^- - EX^-),
$$

$$
B_{nk} = D(2^k) \sum_{i=1}^{n} (\mathbf{1}_{\{X_i^+ \leq D(2^k)\}} - F(D(2^k))),
$$

$$
C_{nk} = \sum_{i=1}^{n} (X_i^+ \mathbf{1}_{\{X_i^+ \leq D(2^k)\}} - \nu_+(D(2^k))).
\tag{3.10}
$$

Then from (3.7) and (3.8),

$$
\begin{aligned}
R_n &\leq |A_n| + \max_{1 \leq j \leq n} |A_j| + D(2^k) \left( |B_{nk}| + \max_{1 \leq j \leq n} |B_{jk}| \right) \\
&\quad + |C_{nk}| + \max_{1 \leq j \leq n} |C_{jk}| + n \int_{D(2^k)}^\infty \bar{F}(y)\,dy,
\end{aligned}
$$

so

$$(3.11) \qquad \max_{1 \le n \le 2^k} R_n \le 2 \max_{1 \le n \le 2^k} |A_n| + 2D(2^k) \max_{1 \le n \le 2^k} |B_{nk}| \\ + 2 \max_{1 \le n \le 2^k} |C_{nk}| + 2^k \int_{D(2^k)}^{\infty} \bar{F}(y)\, dy.$$

The last term on the right-hand side of (3.11) is, by (3.9) and the definition of $W(x)$, not larger than $\delta 2^{\kappa k}$. We will show that the other terms on the right-hand side of (3.11) are $o(2^{\kappa k})$ a.s., as $k \to \infty$.

We need some properties of $D(x)$. Differentiation using the implicit function theorem gives

$$(3.12) \qquad D'(x) = \frac{(1 - \kappa)\delta D^2(x)}{x^{2-\kappa} \int_{[0,D(x)]} y^2 F(dy)}, \qquad x > x_0.$$

Hence, $D(\cdot)$ is strictly increasing and so has a unique increasing inverse $D^{\leftarrow}(x)$ satisfying, for large $x$, $x \ge x_1$, say,

$$(3.13) \qquad D^{\leftarrow}(x) = \left( \frac{\delta x}{W(x)} \right)^{1/(1-\kappa)}.$$

Our next step is to show, under our assumption $I_\kappa(\lambda) < \infty$, that

$$(3.14) \qquad \lim_{x \to \infty} x^{1-2\kappa} W(D(x)) = 0.$$

To see this, write

$$I_\kappa(\lambda) = \int_1^\infty e^{-\lambda y^q / h(y)}\, dy/y,$$

where $q = (2\kappa - 1)/(1 - \kappa) > 0$ and $h(x) = (W(x))^{\kappa/(1-\kappa)}$ is an increasing function. [In fact, differentiation shows that $W(x)$ is increasing and concave.] Now $I_\kappa(\lambda) < \infty$ implies

$$\log 2 \sum_{n \ge 1} e^{-\lambda 2^{(n+1)q}/h(2^n)} \le \sum_{n \ge 1} \int_{2^n}^{2^{n+1}} e^{-\lambda y^q/h(y)}\, dy/y < \infty,$$

thus, $\lim_{n\to\infty} h(2^n)/2^{(n+1)q} = 0$ and so $\lim_{n\to\infty} h(2^n)/2^{(n-1)q} = 0$. Given $x > 0$, choose $n(x)$ so that $2^{n-1} \le x < 2^n$. Then

$$\frac{h(x)}{x^q} \le \frac{h(2^n)}{2^{(n-1)q}} \to 0 \qquad \text{as } x \to \infty,$$

thus,

$$\lim_{x \to \infty} \frac{(W(x))^{\kappa/(1-\kappa)}}{x^{(2\kappa-1)/(1-\kappa)}} = 0.$$

Substituting $x = (D^{\leftarrow}(x))^{1-\kappa}W(x)/\delta$ from (3.13) gives

$$\lim_{x\to\infty} \frac{\delta^{(2\kappa-1)/(1-\kappa)}W(x)}{(D^{\leftarrow}(x))^{2\kappa-1}} = 0,$$

or, equivalently, since $\lim_{x\to\infty}D(x) = \infty$, (3.14) holds, as required.

Now consider first the $C_{nk}$ term in (3.11). By (3.10), $C_{nk}$ is, for each $k$ and $n \leq 2^k$, the sum of $n$ i.i.d. mean 0 r.v.'s, and we can calculate

$$(3.15) \qquad \begin{aligned} \mathrm{Var}(C_{nk}) &= n\,\mathrm{Var}(X_i^+ \mathbf{1}_{\{X_i^+ \leq D(2^k)\}}) \\ &\leq 2^k \int_{[0,D(2^k)]} y^2 F(dy) \leq 2^k W(D(2^k)) = o(2^{2\kappa k}), \end{aligned}$$

where the last estimate follows from (3.14). The inequality $|\mathrm{median}(Y)| \leq \sqrt{2\,\mathrm{Var}\,Y}$ is valid for any mean zero r.v., so we have from (3.15)

$$\max_{1\leq n\leq 2^k}\left|\mathrm{median}\left(\sum_{i=n}^{2^k}(X_i^+ \mathbf{1}_{\{X_i^+\leq D(2^k)\}} - \nu_+(D(2^k)))\right)\right| = o(2^{\kappa k}).$$

Thus, by a version of Lévy's inequality (e.g., [5], page 71), for large enough $k$,

$$(3.16) \qquad P\left(\max_{1\leq n\leq 2^k}|C_{nk}| > 2\delta 2^{\kappa k}\right) \leq 2P(|C_{2^k k}| > \delta 2^{\kappa k}).$$

The summands of $C_{nk}$ are bounded by $2D(2^k)$, so Bernstein's inequality ([5], page 111) gives an upper bound for the last probability as

$$(3.17) \qquad 2\exp\left(\frac{-\delta^2 2^{2\kappa k}}{2(2^k W(D(2^k)) + 2D(2^k)\delta 2^{\kappa k})}\right) = 2\exp\left(\frac{-\delta 2^{\kappa k}}{6D(2^k)}\right),$$

where we used (3.9) to substitute for $W(D(2^k))$. Adding over $k$, we find that

$$\begin{aligned} \sum_{k\geq 1}e^{-\delta 2^{\kappa k}/6D(2^k)} &\leq \sum_{k\geq 1}\int_{2^k}^{2^{k+1}} e^{-\delta y^k/(6\cdot 2^k D(y))}\,dy/y \\ &= \frac{1}{\log 2}\int_2^\infty e^{-\lambda y^\kappa/(\delta^{\kappa/(1-\kappa)}D(y))}\,dy/y, \end{aligned}$$

where in the last we chose $\delta$ so that $\delta = 6\cdot 2^\kappa\lambda/\delta^{\kappa/(1-\kappa)}$, that is, $\delta = (6\lambda 2^\kappa)^{1-\kappa}$. Now change variable to get the last integral as

$$(3.18) \qquad \int_{D(2)}^\infty e^{-\lambda(D^{\leftarrow}(z))^\kappa/(\delta^{\kappa/(1-\kappa)}z)}\frac{dz}{D'(D^{\leftarrow}(z))D^{\leftarrow}(z)}.$$

In view of (3.13), the exponent here is

$$\frac{-\lambda z^{\kappa/(1-\kappa)}}{z(W(z))^{\kappa/(1-\kappa)}} = \frac{-\lambda z^{(2\kappa-1)/(1-\kappa)}}{(W(z))^{\kappa/(1-\kappa)}},$$

as required in (2.6). Also, by (3.12) and (3.13),

$$D'(D^{\leftarrow}(z))D^{\leftarrow}(z) = \frac{(1-\kappa)\delta z^2}{(D^{\leftarrow}(z))^{1-\kappa}\int_{[0,z]}y^2F(dy)}$$

$$= \frac{(1-\kappa)zW(z)}{\int_{[0,z]}y^2F(dy)}$$

$$\geq (1-\kappa)z,$$

where the last follows because $W(z) \geq \int_{[0,z]}y^2F(dy)$; see (2.5). As a result of these two calculations, the integral in (3.18) is bounded by a multiple of $I_\kappa(\lambda)$. Going back to (3.16), we thus have, by the Borel–Cantelli lemma,

$$(3.19) \qquad \limsup_{k\to\infty}\left(\frac{\max_{1\leq n\leq 2^k}|C_{nk}|}{2^{\kappa k}}\right) \leq 2\delta = 2(6\lambda 2^\kappa)^{1-\kappa} \qquad \text{a.s.}$$

Next we have to deal with the $(B)$ term in (3.11). For each $k \geq 1$ and $1 \leq n \leq 2^k$, $B_{nk}/D(2^k)$ is a sum of i.i.d. mean zero r.v.'s bounded by 2 [see (3.10)], and we can calculate

$$\text{Var}(B_{nk}/D(2^k)) = \sum_{i=1}^{n}F(D(2^k))\bar{F}(D(2^k))$$

$$\leq 2^k\bar{F}(D(2^k)) \leq 2^kW(D(2^k))/D^2(2^k)$$

$$= \delta 2^{\kappa k}/D(2^k),$$

using (3.9) for the last equality. Thus, by a similar argument as for the $(C)$ term, involving Lévy's and Bernstein's inequalities,

$$P\left(\max_{1\leq n\leq 2^k}|B_{nk}| > 2\delta 2^{\kappa k}\right) \leq 2P(|B_{2^k k}|/D(2^k) > \delta 2^{\kappa k}/D(2^k))$$

$$\leq 2\exp\left(\frac{-\delta^2 2^{2\kappa k}/D^2(2^k)}{2(\delta 2^{\kappa k}/D(2^k) + 2\delta 2^{\kappa k}/D(2^k))}\right)$$

$$= 2\exp\left(\frac{-\delta 2^{\kappa k}}{6D(2^k)}\right).$$

This is the same bound as in (3.17) and the same argument leading to (3.19) which gives

$$(3.20) \qquad \limsup_{k\to\infty}\left(\frac{\max_{1\leq n\leq 2^k}|B_{nk}|}{2^{\kappa k}}\right) \leq 2(6\lambda 2^\kappa)^{1-\kappa} \qquad \text{a.s.}$$

Finally, for the $(A)$ term in (3.11), we simply use the Marcinkiewicz–Zygmund law to get $A_n = o(n^\kappa)$ a.s., since $E(X^-)^{1/\kappa} < \infty$. So

$$(3.21) \qquad \lim_{k\to\infty}\left(\frac{\max_{1\leq n\leq 2^k}|A_n|}{2^{\kappa k}}\right) = 0 \qquad \text{a.s.}$$

Putting (3.19)–(3.21) into (3.11) gives

$$\limsup_{k \to \infty} \left( \frac{\max_{1 \le n \le 2^k} R_n}{2^{\kappa k}} \right) \le 8(6\lambda 2^{\kappa})^{1-\kappa} + \delta = 9(6\lambda 2^{\kappa})^{1-\kappa} \qquad \text{a.s.}$$

If $m$ is large, choose $k(m)$ so that $2^{k-1} \le m < 2^k$. Then (3.6) follows from

$$\frac{R_m}{m^{\kappa}} \le 2^{\kappa} \frac{\max_{1 \le n \le 2^k} R_n}{2^{\kappa k}} \le 2^{\kappa} 9(6\lambda 2^{\kappa})^{1-\kappa} + o(1) \qquad \text{a.s.,}$$

after letting $\lambda \downarrow \lambda_{\kappa}^*$.  $\square$

Finally we can complete the proof of the forward direction in (2.3). Recall that we are in the case $1/2 < \kappa < 1$, $E|X| < \infty$ and $EX = 0$, and have assumed that $\tau_{\kappa}(r) < \infty$ a.s. for all $r > 0$. Thus, by Lemma 3.1, $\limsup_n R_n/n^{\kappa} = \infty$ a.s. If $E(X^-)^{1/\kappa} < \infty$ and $E(X^+)^{1/\kappa} < \infty$, that is, $E|X|^{1/\kappa} < \infty$, we get $\lim_{n \to \infty} R_n/n^{\kappa} = 0$ a.s. from the Marcinkiewicz–Zygmund law, so we must have $E(X^-)^{1/\kappa} = \infty$ or $E(X^-)^{1/\kappa} < \infty = E(X^+)^{1/\kappa}$. In the latter case we must further have $\lambda_{\kappa}^* = \infty$ by Lemma 3.4. But then $\liminf_{n \to \infty} S_n/n^{\kappa} = -\infty$ a.s. by part (e) of Proposition 2.1.

For the converse part of Theorem 2.1(b), note first that, by its definition, for $r > 0$, $\kappa > 0$, $n = 1, 2, \ldots,$

$$\{\tau_{\kappa}(r) > n\} = \left\{ \max_{0 \le k \le j} S_k - S_j \le r j^{\kappa}, \ 1 \le j \le n \right\}$$
$$\subseteq \{-X_j \le r j^{\kappa}, 1 \le j \le n\},$$

the last following just by taking the term for $k = j - 1$ from the max. So

$$P(\tau_{\kappa}(r) > n) \le \prod_{j=1}^{n} P(X_1 > -r j^{\kappa}) \le \exp\left( -\sum_{j=1}^{n} P(X_1 \le -r j^{\kappa}) \right).$$

Thus, if $\sum_{j \ge 1} P(X_1 \le -r j^{\kappa}) = \infty$, or, equivalently, $E(X^-)^{1/\kappa} = \infty$, then $P(\tau_{\kappa}(r) < \infty) = 1$ for each $r > 0$.

Next, the second condition in (2.3) implies $\limsup_{n \to \infty} R_n/n^{\kappa} = \infty$ a.s., hence, it also implies $P(\tau_{\kappa}(r) < \infty) = 1$ for each $r > 0$ by Lemma 3.1.

This completes the proof of Theorem 2.1.  $\square$

PROOF OF THEOREM 2.2.   (a) For the square root boundary, assume $EX^2 < \infty$ and $EX = 0$.

(i) Introduce the function

$$\phi(x) = 2\left\{ \int_x^{\infty} y\bar{F}(y)\,dy - x \int_x^{\infty} \bar{F}(y)\,dy \right\}$$
$$= 2\int_0^{\infty} y\bar{F}(y+x)\,dy,$$

and define

$$(3.22) \qquad Z_n = R_n^2 - n\sigma^2 + \sum_1^n \phi(R_{i-1}), \qquad n = 1, 2, \ldots, Z_0 = 0.$$

Now, whenever $E|X| < \infty$, we can use (1.4) to write

$$(3.23) \qquad \begin{aligned} E(\Delta_i \mid \mathcal{F}_{i-1}) &= -\int_{(-\infty, R_{i-1}]} y \, dF(y) - R_{i-1} \bar{F}(R_{i-1}) \\ &= -EX + \int_{R_{i-1}}^{\infty} \bar{F}(y) \, dy, \end{aligned}$$

where $\mathcal{F}_i = \sigma(X_1, X_2, \ldots, X_i)$ is the $\sigma$-field generated by $X_1, X_2, \ldots, X_i$, with $\mathcal{F}_0$ as the trivial $\sigma$-field. Using (3.23), and similarly calculating

$$\begin{aligned} E(\Delta_i^2 \mid \mathcal{F}_{i-1}) &= \int_{(-\infty, R_{i-1}]} y^2 \, dF(y) + R_{i-1}^2 \bar{F}(R_{i-1}) \\ &= \sigma^2 - 2\int_{R_{i-1}}^{\infty} y \bar{F}(y) \, dy, \end{aligned}$$

we can write

$$Z_n = R_n^2 - \sum_1^n E(\Delta_i^2 \mid \mathcal{F}_{i-1}) - 2\sum_1^n R_{i-1} E(\Delta_i \mid \mathcal{F}_{i-1}).$$

From this, it is easy to check that $Z$ is a martingale. Now fix $r > 0$ and $m > 0$ and write $\tau$ for $\tau_{1/2}(\sigma r)$ and $\tau^m = m \wedge \tau$. This is a bounded stopping time, so by Doob's theorem (e.g., [5]), $EZ_{\tau^m} = 0$, and thus, from (3.22),

$$(3.24) \qquad \sigma^2 E\tau^m = ER_{\tau^m}^2 + E\sum_1^{\tau^m} \phi(R_{i-1}) \geq ER_{\tau^m}^2 + \phi(0).$$

Suppose now that $E\tau < \infty$. By monotone convergence, $E\tau^m \to E\tau$ as $m \to \infty$, while

$$\liminf_{m\to\infty} ER_{\tau^m}^2 \geq ER_\tau^2 \geq r^2 \sigma^2 E\tau$$

by Fatou's lemma. Thus, we can let $m \to \infty$ in (3.24) to get $\sigma^2(1 - r^2)E\tau \geq \phi(0) > 0$. This is impossible if $r \geq 1$, so in this case we must have $E\tau = \infty$.

(ii) We now take $0 < r < 1$, assume $E\tau = \infty$, and establish a contradiction. Assume the truth of the following statement: for any $\varepsilon > 0$, there is an $m_\varepsilon$ such that

$$(3.25) \qquad E\sum_1^{\tau^m} \phi(R_{i-1}) \leq \varepsilon E\tau^m \qquad \text{for all } m \geq m_\varepsilon.$$

Note that $R_{\tau^m} = R_{\tau^m-1} + \Delta_{\tau^m} \leq \sigma r \sqrt{\tau^m} + \Delta_{\tau^m}$, and choose $\varepsilon \in (0, \sigma^2)$. Then for any $m \geq m_\varepsilon$, we have, using the equality in (3.24), and (3.25),

$$
\begin{aligned}
\sigma^2 E \tau^m &= E R_{\tau^m}^2 + E \sum_1^{\tau^m} \phi(R_{i-1}) \\
&\leq \sigma^2 r^2 E \tau^m + E \Delta_{\tau^m}^2 + 2\sigma r E(\sqrt{\tau^m} \Delta_{\tau^m}) + \varepsilon E \tau^m \\
&\leq \sigma^2 r^2 E \tau^m + E \Delta_{\tau^m}^2 + 2\sigma r \sqrt{E \tau^m} \sqrt{E \Delta_{\tau^m}^2} + \varepsilon E \tau^m.
\end{aligned}
$$

Thus,

$$
(3.26) \qquad (\sigma^2 - \varepsilon) E \tau^m \leq (\sqrt{E \Delta_{\tau^m}^2} + \sigma r \sqrt{E \tau^m})^2.
$$

From this, we see that the ratio $E \Delta_{\tau^m}^2 / E \tau^m$ is bounded below when $\varepsilon$ is small enough, $m \geq m_\varepsilon$, and $r < 1$. The contradiction will follow by showing that $E \Delta_{\tau^m}^2 / E \tau^m \to 0$ as $m \to \infty$.

To see this, take any $\delta > 0$. First note that we can choose $M = M(\varepsilon, \delta) \geq m_\varepsilon$ so large that, whenever $m \geq M$,

$$
(3.27) \qquad \max_{i \geq 1} E(\Delta_i^2 \mathbf{1}_{\{\Delta_i^2 > \varepsilon E \tau^m\}} \mid \mathcal{F}_{i-1}) \leq \delta \qquad \text{a.s.}
$$

This can be demonstrated as follows. Since $EX^2 < \infty$, given $\delta > 0$, we can choose $y_0(\delta)$ so large that

$$
(3.28) \qquad \int_{|z|>y} z^2 \, dF(z) + \sup_{z>y} z^2 \bar{F}(z) \leq \delta \qquad \text{for all } y \geq y_0.
$$

Since we assumed $E\tau = \infty$, we have $\lim_{m\to\infty} E\tau^m = \infty$. So we can also choose $M(\varepsilon, \delta)$ so large that $\sqrt{\varepsilon E \tau^m} \geq y_0$ when $m \geq M$. Now for any $a > 0$, using the representation (1.4),

$$
\mathbf{1}_{\{\Delta_i^2 > a\}} = \mathbf{1}_{\{X_i^2 > a\}} \mathbf{1}_{\{X_i \leq R_{i-1}\}} + \mathbf{1}_{\{R_{i-1}^2 > a\}} \mathbf{1}_{\{X_i > R_{i-1}\}},
$$

so

$$
\Delta_i \mathbf{1}_{\{\Delta_i^2 > a\}} = -X_i \mathbf{1}_{\{X_i^2 > a\}} \mathbf{1}_{\{X_i \leq R_{i-1}\}} - R_{i-1} \mathbf{1}_{\{R_{i-1}^2 > a\}} \mathbf{1}_{\{X_i > R_{i-1}\}}.
$$

Hence,

$$
\begin{aligned}
(3.29) \quad & E(\Delta_i^2 \mathbf{1}_{\{\Delta_i^2 > a\}} \mid \mathcal{F}_{i-1}) \\
&= \int_{|y|>\sqrt{a}} y^2 \mathbf{1}_{\{y \leq R_{i-1}\}} \, dF(y) + R_{i-1}^2 \mathbf{1}_{\{R_{i-1}>\sqrt{a}\}} \bar{F}(R_{i-1}) \\
&\leq \int_{|y|>\sqrt{a}} y^2 \, dF(y) + \sup_{y>\sqrt{a}} y^2 \bar{F}(y).
\end{aligned}
$$

Substituting $a = \varepsilon E \tau^m$, we have $\sqrt{a} \geq y_0$ when $m \geq M$, so (3.29) gives (3.27) via (3.28), when $m \geq M$. From (3.27), still with $a = \varepsilon E \tau^m$, we deduce

$$E \Delta_{\tau^m}^2 \mathbf{1}_{\{\Delta_{\tau^m}^2 > a\}} \leq E \sum_{i=1}^{\tau^m} \Delta_i^2 \mathbf{1}_{\{\Delta_i^2 > a\}}$$

$$= E \sum_{i \geq 1} E(\mathbf{1}_{\{\tau^m > i-1\}} \Delta_i^2 \mathbf{1}_{\{\Delta_i^2 > a\}} \mid \mathcal{F}_{i-1})$$

$$= E \sum_{i=1}^{\tau^m} E(\Delta_i^2 \mathbf{1}_{\{\Delta_i^2 > a\}} | \mathcal{F}_{i-1})$$

$$\leq \delta E \tau^m,$$

for $m \geq M$. We also have $E \Delta_{\tau^m}^2 \mathbf{1}_{\{\Delta_{\tau^m}^2 \leq a\}} \leq a = \varepsilon E \tau^m$. (3.26) then gives

$$(\sigma^2 - \varepsilon) E \tau^m \leq (\sqrt{\varepsilon + \delta} + \sigma r)^2 E \tau^m$$

for $m \geq M$, which is impossible for $\varepsilon$ and $\delta$ small enough, when $r < 1$. So to complete the proof, it suffices to prove (3.25).

Note first that $\phi(x)/2 \leq \sigma_+^2 := E(X^+)^2$ for all $x \geq 0$, and $\phi(x) \downarrow 0$ as $x \to \infty$. Fix $\varepsilon > 0$ and choose $K_\varepsilon < \infty$ such that $\phi(K_\varepsilon) \leq \varepsilon/3$. Then we have the bound

$$\sum_1^n \phi(R_{i-1}) \leq \tfrac{1}{3} n \varepsilon + 2 \sigma_+^2 \sum_1^n \mathbf{1}_{\{R_{i-1} \leq K_\varepsilon\}}.$$

Define

$$N^{(\varepsilon)} = \max \left( n : \sum_1^n \mathbf{1}_{\{R_{i-1} \leq K_\varepsilon\}} \geq \frac{n\varepsilon}{6\sigma_+^2} \right).$$

Then it suffices to show that $E N^{(\varepsilon)} < \infty$, since this gives

$$E \sum_1^{\tau^m} \phi(R_{i-1}) \mathbf{1}_{\{\tau^m \leq N_\varepsilon\}} \leq 2 \sigma_+^2 E N_\varepsilon = o(E \tau^m),$$

while

$$E \sum_1^{\tau^m} \phi(R_{i-1}) \mathbf{1}_{\{\tau^m > N_\varepsilon\}} \leq \frac{\varepsilon}{3} E \tau^m + 2 \sigma_+^2 E \left( \frac{\tau^m \varepsilon}{6 \sigma_+^2} \right) = \frac{2\varepsilon}{3} E \tau^m.$$

To show that $E N^{(\varepsilon)} < \infty$, introduce the r.v.s $\alpha_n$, $\beta_n$, $n \geq 1$, given recursively by

$$\alpha_1 = \min\{n \geq 1 : R_n > K_\varepsilon\},$$

$$\beta_1 = \min\{n \geq 1 : R_{\alpha_1 + n} \leq K_\varepsilon\},$$

$$\gamma_1 = \alpha_1 + \beta_1,$$

and, for $i = 2, 3, \ldots,$

$$\alpha_i = \min\{n \geq 1 : R_{\gamma_{i-1}+n} > K_\varepsilon\},$$

$$\beta_i = \min\{n \geq 1 : R_{\gamma_{i-1}+\alpha_i+n} \leq K_\varepsilon\},$$

$$\gamma_i = \gamma_{i-1} + \alpha_i + \beta_i.$$

In view of Proposition 1.1, the $\alpha_i$ and $\beta_i$ are finite, a.s. Then, by construction,

$$(3.30) \qquad \sum_1^n \mathbf{1}_{\{R_{i-1} \leq K_\varepsilon\}} = u_n := \sum_1^{d_n} \alpha_i + (n - \gamma_{d_n}),$$

where $d_n = \max\{k : \gamma_k \leq n\}$. Now write $\widetilde{\varepsilon} = \varepsilon\sigma^2/(4\sigma_+^2)$, assume without loss of generality that $\widetilde{\varepsilon} < 1$, and note that the maximum values of $n^{-1}u_n$ occur when $n = \gamma_k + \alpha_{k+1}$ for some $k \geq 0$, that is, when $\gamma_{d_n} = n - \alpha_{k+1}$, at which times $u_n$ has the value $\sum_{i=1}^{k+1} \alpha_i$. So

$$N^{(\varepsilon)} = \max\{n : u_n \geq n\widetilde{\varepsilon}\}$$

$$\leq \max\left\{\gamma_k + \alpha_{k+1} : \sum_{i=1}^{k+1} \alpha_i \geq \widetilde{\varepsilon}(\gamma_k + \alpha_{k+1})\right\}$$

$$\leq \frac{1}{\widetilde{\varepsilon}} \max\left\{\sum_{i=1}^{k+1} \alpha_i : \sum_{i=1}^{k+1} \alpha_i \geq \widetilde{\varepsilon}\left(\sum_{i=1}^{k+1} \alpha_i + \sum_{i=1}^{k} \beta_i\right)\right\}$$

$$= \frac{1}{\widetilde{\varepsilon}} \max\left\{\sum_{i=1}^{k+1} \alpha_i : \sum_{i=1}^{k} (1-\widetilde{\varepsilon})\alpha_i - \widetilde{\varepsilon}\beta_i \geq (\widetilde{\varepsilon}-1)\alpha_{k+1}\right\}.$$

Thus, writing $Y_i = \widetilde{\varepsilon}\beta_i - (1-\widetilde{\varepsilon})\alpha_i$ and $k^* = \max\{k : \sum_1^k Y_i \leq (1-\widetilde{\varepsilon})\alpha_{k+1}\}$, we have, for any $c > 0$,

$$(3.31) \qquad \begin{aligned} P(N^{(\varepsilon)} \geq m/\widetilde{\varepsilon}) &\leq P\left(\sum_1^{k^*+1} \alpha_i \geq m\right) \\ &\leq P(k^* \geq mc - 1) + P\left(\sum_1^{mc} \alpha_i \geq m\right). \end{aligned}$$

Now $\sum_1^{mc} \alpha_i \leq \sum_1^{mc} \widetilde{\alpha}_i$, where $\widetilde{\alpha}_1, \widetilde{\alpha}_2, \ldots$ are i.i.d. with the distribution of the time that $R$, starting from 0, crosses the level $K_\varepsilon$. Part (a) of Theorem 2.1 shows that $Ee^{\lambda\widetilde{\alpha}_1} < \infty$ for some $\lambda > 0$, so using a standard exponential bound and choosing $c < \lambda/\log Ee^{-\lambda\widetilde{\alpha}_1}$, we see that the second term in (3.31) is summable. On the other hand, we have

$$(3.32) \quad \beta_k \geq \widetilde{\beta}_k := \min\{n : R_{\gamma_{k-1}+\alpha_k+n} \leq R_{\gamma_{k-1}+\alpha_k}\} \geq \min\{n : \hat{S}_n \geq 0\},$$

where $\hat{S}_n = S_{\gamma_{k-1}+\alpha_k+n} - S_{\gamma_{k-1}+\alpha_k}$, $n \geq 0$, and the $\widetilde{\beta}_n$ are an i.i.d. sequence with infinite mean since $\liminf_n S_n = -\infty$ a.s.; see Theorem II.9.1(iii) of [14], page 66. Thus, $\widetilde{Y}_i := \widetilde{\varepsilon}\widetilde{\beta}_i - (1-\widetilde{\varepsilon})\alpha_i$ are the i.i.d. steps of a random walk that drifts to $+\infty$ a.s.; see, for example, Theorem II.8.3(i) in [14], page 64.

Then with $A(y) := E((Y_1 \wedge y) \vee (-y))$, write

$$\sum_{j \geq 1} P(k^* \geq j) = \sum_{j \geq 1} P\left(\text{for some } k \geq j, \sum_{i=1}^{k} \widetilde{Y}_i \leq (1-\widetilde{\varepsilon})\alpha_{k+1}\right)$$

$$\leq \sum_{j \geq 1}\sum_{k \geq j}\sum_{a \geq 1} P\left(\sum_{i=1}^{k} \widetilde{Y}_i \leq (1-\widetilde{\varepsilon})a\right)P(\alpha_{k+1}=a)$$

$$= \sum_{a \geq 1}\sum_{k \geq 1} k P\left(\sum_{i=1}^{k} \widetilde{Y}_i \leq (1-\widetilde{\varepsilon})a\right)P(\alpha_{k+1}=a)$$

$$\leq c_1 + c_2 \sum_{a \geq a_0}\left(\frac{(1-\widetilde{\varepsilon})a}{A((1-\widetilde{\varepsilon})a)}\right)^2 P(\alpha_{k+1}=a)$$

$$\leq c_1 + c_3 E\alpha_1^2 < \infty.$$

Here the $c_i$ and $a_0$ are positive constants, $A(y)$ is bounded away from 0 for $a \geq a_0$, say [in fact, $\lim_{y \to \infty} A(y) = \infty$ since $E\widetilde{Y}_1 = +\infty$], and the inequality in the fourth line follows from Theorem 2.2 of [20]. Thus, $k^*$ has finite mean, and the result follows from (3.31).

(b) For the linear case, assume $E|X| < \infty$ and $EX < 0$.

(i) We first show that $E\tau_1(r) < \infty$ for $r < -EX$. We have

$$\tau_1(r) = \min\{n \geq 1 : R_n > rn\} \leq \min\{n \geq 1 : S_n < -rn\}$$

$$= \min\left\{n \geq 1 : \sum_{i=1}^{n}(X_i - EX - (|EX|-r)) < 0\right\},$$

so $\tau_1(r)$ does not exceed the first strict decreasing ladder time of a random walk which has negative drift. Thus, $E\tau_1(r) < \infty$ in this case.

(ii) Now take $r > |EX|$. If $E\tau_1(r) < \infty$, then $\tau_1(r) < \infty$ a.s., so $\limsup_n R_n / n > r$ a.s. by Lemma 3.1, contradicting $\lim_{n \to \infty} R_n / n = |EX| < r$ a.s., which follows from the strong law of large numbers for $S_n$.

(c) Assume $E|X| < \infty$ and $EX < 0$, and in addition, that $E(X^+)^2 < \infty$. We will show that $E\tau_1(r) = \infty$ when $r = |EX|$. This follows immediately from the next lemma, which proves a little more.

LEMMA 3.5.   *Let $S$ be a random walk with steps $X$ having $E|X| < \infty$, $EX = \mu < 0$ and $E(X^+)^2 < \infty$, and for the corresponding reflected process*

$R_n = \max_{0 \le i \le n} S_i - S_n$, $R_0 = 0$, *write*

(3.33)                    $T_a = \min(n \ge 1 : R_n > n|\mu| + a).$

*Then* $ET_a = \infty$ *for* $a \ge 0$.

REMARK.  Of course, $ET_0 = \infty$ implies $ET_a = \infty$ for $a \ge 0$, but it does not seem possible to prove it without considering the case $a > 0$.

PROOF OF LEMMA 3.5.  Since $EX = \mu < 0$, but $P(X \ge 0) > 0$, there is probability mass above $\mu$, so we can assume $P(X \ge -|\mu| + \delta) > c > 0$ for some $\delta \in (0, |\mu|)$. First we show the required result holds for sufficiently large $a$. Note that $R_n = S_n^* + n|\mu| - \widetilde{S}_n$, where $\widetilde{S}$ is a zero-mean random walk and $S_n^* = \max_{0 \le i \le n} S_i \le S_\infty^*$, where $b := ES_\infty^* < \infty$ since $E(X^+)^2 < \infty$; see [19]. Assuming $ET_a < \infty$, we get

$$0 = E\widetilde{S}_{T_a} = ES_{T_a}^* + |\mu|ET_a - ER_{T_a}$$
$$\le ES_\infty^* + |\mu|ET_a - (a + |\mu|ET_a) = b - a.$$

This is a contradiction when $a > b$, so $ET_a = \infty$ for $a > b$.

Next, observe that $ET_a \ge P(X \ge 0)ET_{a+|\mu|} = cET_{a+|\mu|}$, for a $c > 0$. Thus, if $ET_0$ were finite, $ET_{n|\mu|}$ would also be finite for $n = 1, 2, \ldots$. This proves the lemma.  □

With this, the proof of Theorem 2.2 is complete.  □

**Acknowledgments.**  We are grateful to Chris Wetherell for very competent LaTex typing. Thanks also to two referees for an extremely careful reading of the paper, many helpful suggestions and for pointing out some errors and obscurities; and likewise also to Mr. M. Savov.

School of Mathematics
University of Manchester
Manchester M60 1QD
UK
E-mail: rad@maths.man.ac.uk

Centre for Mathematical Analysis
  and School of Finance and Applied Statistics
Australian National University
Canberra, ACT
Australia
E-mail: Ross.Maller@anu.edu.au